\documentclass[11pt]{amsart}
\usepackage[all]{xy}
\usepackage{amsmath}
\usepackage{amsthm}
\usepackage{amsfonts}
\usepackage{amssymb}
\usepackage{amscd}
\usepackage{epic}
\usepackage{eepic}
\usepackage{verbatim}
\newtheorem{theorem}{Theorem}[section]
\newtheorem*{theorem*}{Theorem}
\newtheorem{proposition}[theorem]{Proposition}

 \theoremstyle{definition}

\newtheorem{dfn}[theorem]{Definition}

\newcommand{\ind}{\operatorname{ind}}

\newcommand{\gPGL}{\operatorname{\mathbf{PGL}}}

\newcommand{\Z}{\mathbb{Z}}

\newcommand{\R}{\mathbb{R}}

\newcommand{\CH}{\operatorname{CH}}

\newcommand{\ra}{\rightarrow}

\newcommand{\cdim}{\operatorname{cdim}}
\newcommand{\al}{\alpha}
\newcommand{\ds}{\oplus}
\newcommand{\mult}{\operatorname{mult}}
\newcommand{\de}{\delta}
\newcommand{\lt}{\leadsto}
\newcommand{\im}{\operatorname{im}}
\newcommand{\trans}{\operatorname{transpose}}
\newcommand{\Spec}{\operatorname{Spec}}

\newcommand{\rat}{\dashrightarrow}
\renewcommand{\phi}{\varphi}

\newcommand{\ten}{\otimes}
\newcommand{\cha}{\operatorname{char}}

\title{Canonical dimension of projective $\gPGL_1(A)$-homogeneous varieties}
\author{Bryant Mathews}
\date{}
\address{Department of Mathematics, University of California,
        Los Angeles, CA 90095-1555} \email {bmathews@math.ucla.edu}

\begin{document}
\begin{abstract}
Let $A$ be a central division algebra over a field $F$ with $\ind A =n$.  For integers $1\leq d_1<d_2<\cdots < d_k\leq n-1$, let $X_{d_1,d_2,\ldots,d_k}(A)$ be the variety of flags of right ideals $I_1\subset I_2\subset\cdots\subset I_k$ of $A$ with $I_i$ of reduced dimension $d_i$.  In computing canonical $p$-dimension of such varieties, for $p$ prime, we can reduce to the case of generalized Severi-Brauer varieties $X_e(A)$ with $\ind A$ a power of $p$ divisible by $e$.  We prove that canonical $2$-dimension (and hence canonical dimension) equals dimension for all $X_e(A)$ with $\ind A=2e$ a power of $2$.
\end{abstract}
\maketitle

\section{Canonical $p$-dimension}

We begin by recalling the definitions of canonical $p$-dimension, $p$-incompressibility, and equivalence.

Let $X$ be a scheme over a field $F$, and let $p$ be a prime or zero.  A field extension $K$ of $F$ is called a \emph{splitting field of $X$} (or is said to \emph{split $X$}) if $X(K)\neq\emptyset$.  A splitting field $K$ is called \emph{$p$-generic} if, for any splitting field $L$ of $X$, there is an $F$-place $K\rightharpoonup L'$ for some finite extension $L'/L$ of degree prime to $p$.  In particular, $K$ is $0$-generic if for any splitting field $L$ there is an $F$-place $K\rightharpoonup L$.

The canonical $p$-dimension of a scheme $X$ over $F$ was originally defined \cite{BR05,KM06} as the minimal transcendence degree of a $p$-generic splitting field $K$ of $X$.  When $X$ is a smooth complete variety, the original algebraic definition is equivalent to the following geometric one \cite{KM06,M}.
\begin{dfn}\label{cddef}
Let $X$ be a smooth complete variety over $F$.  The \emph{canonical $p$-dimension $\cdim_p(X)$ of $X$} is the minimal dimension of the image of a morphism $X'\ra X$, where $X'$ is a variety over $F$ admitting a dominant morphism $X'\ra X$ with $F(X')/F(X)$ finite of degree prime to $p$.  The canonical $0$-dimension of $X$ is thus the minimal dimension of the image of a rational morphism $X\rat X$.
\end{dfn}
In the case $p=0$, we will drop the $p$ and speak simply of \emph{generic} splitting fields and canonical \emph{dimension} $\cdim(X)$.

For a third definition of canonical $p$-dimension as the essential $p$-dimension of the detection functor of a scheme $X$, we refer the reader to Merkurjev's comprehensive exposition \cite{M} of essential dimension. \\

For a smooth complete variety $X$, the inequalities $$\cdim_p(X)\leq \cdim (X)\leq \dim (X)$$ are clear from Definition \ref{cddef}.  Note also that if $X$ has a rational point, then $\cdim(X)=0$ (though the converse is not true).

\begin{dfn}
When a smooth complete variety $X$ has canonical $p$-dimension as large as possible, namely $\cdim_p(X)=\dim(X)$, we say that $X$ is \emph{$p$-incompressible}.
\end{dfn}

It follows immediately that if $X$ is $p$-incompressible, it is also \emph{incompressible} (i.e. $0$-incompressible). \\

When two schemes $X$ and $Y$ over a field $F$ have the same class of splitting fields, we call them equivalent and write $X\sim Y$.  In this case $$\cdim_p(X)=\cdim_p(Y)$$ for all $p$.  If $X$ and $Y$ are smooth complete varieties, then they are equivalent if and only if there exist rational maps $X\rat Y$ and $Y\rat X$.

\section{Reductions}

Let $A$ be a central division algebra over a field $F$ with $\ind A =n$.  We consider the problem of computing the canonical $p$-dimension of the following varieties.

\begin{dfn}
For integers $1\leq d_1<d_2<\cdots < d_k\leq n-1$, define $X_{d_1,d_2,\ldots,d_k}(A)$ to be the variety of flags of right ideals $I_1\subset I_2\subset\cdots\subset I_k$ of $A$ with $I_i$ of reduced dimension $d_i$.  When the algebra $A$ is understood, we write simply $X_{d_1,d_2,\ldots,d_k}$.
\end{dfn}

When $k=1$ we get the generalized Severi-Brauer varieties $X_d(A)$ of $A$.  In particular, $X_1(A)$ is the Severi-Brauer variety of $A$.

It is known \cite[Th. 1.17]{KMRT98} that the generalized Severi-Brauer variety $X_{d_1}(A)$ has a rational point over an extension field $K/F$ if and only if the index $\ind A_K$ divides $d_1$.  As a consequence, $X_{d_1}(A)\sim X_{d}(A)$, where $d:=\gcd(\ind A, d_1).$  We record the easy generalization of this fact to varieties $X_{d_1,d_2,\ldots,d_k}(A)$.

\begin{proposition}\label{flagequiv}
If $d:=\gcd(\ind A,d_1,d_2,\ldots,d_k)$, then $$X_{d_1,d_2,\ldots,d_k}(A)\sim X_d(A)$$ and thus, for any $p$,  $$\cdim_p(X_{d_1,d_2,\ldots,d_k}(A))=\cdim_p(X_d(A)).$$
\end{proposition}

\begin{proof}
If $X_{d_1,d_2,\ldots,d_k}(A)$ has a rational point over an extension field $K/F$, then by definition $A_K$ has right ideals of reduced dimensions $d_1,d_2,\ldots,d_k$.  This is the case if and only if $\ind A_K$ divides each of the $d_i$, or equivalently, $\ind A_K$ divides $d$ (since $\ind A_K$ always divides $\ind A$).

Reading the argument backwards, $\ind A_K$ dividing $d$ implies the existence of right ideals $I_1,I_2,\ldots,I_k$ in $A_K$ with reduced dimensions $d_1,d_2,\ldots,d_k$.  In fact, the $I_1,\ldots,I_k$ can be chosen to form a flag.  Suppose $d_i=m_i\ind A_K$ and $A_K\simeq M_t(D)$ for some division algebra $D$.  Then we take $I_i$ to be the set of matrices in $M_t(D)$ whose $t-m_i$ last rows are zero.
\end{proof}

\emph{Hence it is enough to compute $\cdim_p(X_d(A))$ for $d$ dividing $\ind A$.}

If the index of $A$ factors as $\ind A=q_1q_2\cdots q_r$ with the $q_j$ powers of distinct primes $p_j$, then there exist central division algebras $A_j$ of index $q_j$ for $j=1,\ldots,r$ such that $$A\simeq A_1\ten A_2\ten \cdots\ten A_r.$$

\begin{proposition}\label{decomp}
Given a positive integer $1\leq d\leq \ind A-1$, with $q_j$ as above, define $e_j:=\gcd(d,q_j)$ for $j=1,\ldots,r$.  Then $$X_d(A)\sim X_{e_1}(A_1)\times X_{e_2}(A_2)\times\cdots\times X_{e_r}(A_r)$$ and thus, for any $p$, $$\cdim_p(X_d(A)) = \cdim_p(X_{e_1}(A_1)\times X_{e_2}(A_2)\times\cdots\times X_{e_r}(A_r)).$$
\end{proposition}

\begin{proof}
The variety $X_d(A)$ has a rational point over an extension field $K/F$ if and only if $\ind A_K$ divides $d$.  Because $$\ind A_K = (\ind (A_1)_K)\cdots(\ind (A_r)_K),$$ this condition is equivalent to $\ind (A_j)_K$ dividing $d$ for all $j$, or to $\ind (A_j)_K$ dividing $e_j$ for all $j$ (since $\ind (A_j)_K$ always divides $\ind A_j=q_j$).  This holds if and only if each $X_{e_j}(A_j)$ has a rational point over $K$, which is equivalent to the product of the $X_{e_j}(A_j)$ having a rational point over $K$.
\end{proof}

The proposition gives the following upper bound on canonical $p$-dimension:
\begin{equation}\label{bound}
\cdim_p(X_d(A)) \leq \dim\prod_{j=1}^r{X_{e_j}(A_j)} = \sum_{j=1}^r{\dim X_{e_j}(A_j)} = \sum_{j=1}^r{e_j(q_j-e_j)}.
\end{equation}

If $p$ is prime, then there exists a finite, $p$-coprime extension $K$ of $F$ which splits the algebras $A_j$ for all $j$ with $p_j\neq p$.  Since canonical $p$-dimension does not change under such an extension \cite[Prop. 1.5 (2)]{M}, $\cdim_p(X_d(A)) = 0$ unless some $p_s=p$, in which case $$\cdim_p(X_d(A))=\cdim_p(X_{e_s}(A_s)).$$

\emph{We see that it is enough, when $p$ is prime, to compute the canonical $p$-dimension of varieties of the form $X_e(A)$ with $\ind A$ a prime power divisible by $e$. When $p=0$, it is enough to compute the canonical dimension of products of such varieties.}

\section{Known results for Severi-Brauer varieties}

We now recall what is already known about the canonical $p$-dimension of Severi-Brauer varieties $X_1(A)$, the $d=1$ case.

For any $p$, if $d=1$ in (\ref{bound}) above, then all of the $e_j=1$, and the upper bound becomes
\begin{equation}\label{boundd}
\cdim_p(X_1(A)) \leq \sum_{j=1}^r{(q_j-1)}.
\end{equation}
In the special case $r$=1 and $p=p_1$, it is shown in \cite[Th. 11.4]{BR05}, based on Karpenko's \cite[Th. 2.1]{K00}, that the inequality (\ref{boundd}) is actually an equality.  Thus, for general $A$, we have
$$ \cdim_{p_j}(X_1(A))=\cdim_{p_j}(X_1(A_j))=q_j-1$$
for $j=1,2,\ldots,r$, while $\cdim_p(X_1(A))=0$ for all other {\em primes} $p$ \cite[Ex. 5.10]{KM06}.

Now let $p=0$, $d=1$.  When $r=1$, we again have equality in (\ref{boundd}), since canonical dimension is bounded below by canonical $p$-dimension for every prime $p$.  In \cite[Th. 1.3]{CKM07}, (\ref{boundd}) is proven also to be an equality in the case $\ind A = 6$ (i.e. $r=2$, $q_1=2$, $q_2=3$) provided that $\cha F=0$.  The authors of \cite{CKM07} suggest that equality may indeed hold for any $A$ when $p=0$, $d=1$.

\section{$2$-Incompressibility of $X_e(A)$ for $\ind A=2e$ a power of $2$}

If $A$ is a central division algebra with $\ind A=4$, the variety $X_2(A)$ is known to be $2$-incompressible.  Indeed, if the exponent of $A$ is 2, then $X_2(A)$ is isomorphic to a $4$-dimensional projective quadric hypersurface called the \emph{Albert quadric} of $A$ \cite[\S 5.2]{M99}.  Such a quadric has first Witt index $1$ \cite[p. 93]{V}, hence is $2$-incompressible by \cite[Th. 90.2]{EKM08}.  If the exponent of $A$ is 4, we can reduce to the exponent $2$ case by extending to the function field of the Severi-Brauer variety of $A\otimes A$.

In what follows, we show $2$-incompressibility for an infinite family of varieties which includes the varieties of the form $X_2(A)$ (with $\ind A=4$) mentioned above.

\begin{theorem}
Let $e=2^a$, $a\geq 1$.  For a central division algebra $A$ with $\ind A = 2e$, the variety $X_e:=X_e(A)$ is $2$-incompressible.  Thus $$\cdim_2(X_e) = \cdim(X_e) = \dim(X_e) = e(2e-e)= e^2=4^a.$$
\end{theorem}

We briefly recall some terminology from \cite[\S 62 and \S 75]{EKM08}. Let $X$ and $Y$ be schemes with $\dim X=e$.  A \emph{correspondence of degree zero $\de:X\lt Y$ from $X$ to $Y$} is just a cycle $\de\in\CH_e(X\times Y)$.  The \emph{multiplicity} $\mult(\de)$ of such a $\de$ is the integer satisfying $\mult(\de)\cdot [X] = p_*(\de)$, where $p_*$ is the push-forward homomorphism
$$p_*:\CH_e(X\times Y)\ra \CH_e(X) = \mathbb{Z}\cdot [X].$$  The exchange isomorphism $X\times Y\ra Y\times X$ induces an isomorphism
$$\CH_e(X\times Y)\ra \CH_e(Y\times X)$$
sending a cycle $\de$ to its \emph{transpose} $\de^t$.

To prove that a variety $X$ is $2$-incompressible, it suffices to show that for any correspondence $\de:X\lt X$ of degree zero,
\begin{equation}\label{mod}
\mult(\de)\equiv\mult(\de^t) \pmod{2}.
\end{equation}
Indeed, suppose we have $f:X'\ra X$ and a dominant $g:X'\ra X$ with $F(X')/F(X)$ finite of odd degree.  Let $\de\in\CH(X\times X)$ be the pushforward of the class $[X']$ along the induced morphism $(g,f):X'\ra X\times X$.  By assumption, $\mult(\de)$ is odd, so by (\ref{mod}) we have that $\mult(\de^t)$ is odd.  It follows that $f_*([X'])$ is an odd multiple of $[X]$ and in particular is nonzero, so $f$ is dominant.

We will check that the condition (\ref{mod}) holds for the variety $X_e$.  A correspondence of degree zero $\de:X_e\lt X_e$ is just an element of $\CH_{e^2}(X_e\times X_e)$.  Using the method of Chernousov and Merkurjev described in \cite{CM06}, we can decompose the Chow motive of $X_e\times X_e$ as follows.  See also \cite{CSM05} for examples of similar computations.

We first realize $X_e$ as a projective homogeneous variety.  Let $n:=\ind A = 2e=2^{a+1}$.  Let $G$ denote the group $\gPGL_1(A)$, and let $\Pi$ be a set of simple roots for the root system $\Sigma$ of $G$.  If $\varepsilon_1,\ldots,\varepsilon_n$ are the standard basis vectors of $\mathbb{R}^n$, we may take $$\Pi=\{\al_1\!:=\!\varepsilon_1\!-\!\varepsilon_2,\ldots,\al_{n-1}\!:= \!\varepsilon_{n-1}\!-\!\varepsilon_n\}.$$  Then $X_e$ is a projective $G$-homogeneous variety, namely the variety of all parabolic subgroups of $G$ of type $S$, for the subset $S=\Pi\backslash\{\al_e\}$ of the set of simple roots.

Let $W$ denote the Weyl group of the root system $\Sigma$.  There are $e+1$ double cosets $D\in W_P\backslash W/W_P$ with representatives $w$ as follows, where $w_{\al_k}$ denotes the reflection induced by the root $\al_k$.

\begin{gather*}
1 \\
w_{\al_e} \\
(w_{\al_{e}}w_{\al_{e-1}})(w_{\al_{e+1}}w_{\al_e)} \\
(w_{\al_{e}}w_{\al_{e-1}}w_{\al_{e-2}})(w_{\al_{e+1}}w_{\al_{e}}w_{\al_{e-1}})(w_{\al_{e+2}}w_{\al_{e+1}}w_{\al_{e}}) \\
\vdots \\
(w_{\al_{e}}\cdots w_{\al_{1}})\cdots (w_{\al_{2e-1}}\cdots w_{\al_{e}})
\end{gather*}

The subset of $\Pi$ associated to $w=1$ is of course $S=\Pi\backslash\{\al_e\}$.  The general nontrivial representative $$w=w^{-1}=(w_{\al_{e}}\cdots w_{\al_{e-i}})\cdots(w_{\al_{e+i}}\cdots w_{\al_{e}}),$$ for $i\in\{0,\ldots,e-1\}$, has the effect on $\R^n$ of switching the tuple of standard basis vectors $(\varepsilon_{e-i},\ldots,\varepsilon_e)$ with the tuple $(\varepsilon_{e+1},\ldots,\varepsilon_{e+1+i})$.  The resulting subset associated to $w$ is therefore $$\Pi\backslash\{\al_{e-(i+1)},\al_e,\al_{e+(i+1)}\}$$ for $i=0,\ldots,e-2$ and $\Pi\backslash\{\al_e\}$ for $i=e-1$.

From Theorem 6.3 of \cite{CM06}, we deduce the following decomposition of the Chow motive of $X_e\times X_e$, where the relation between the indices $i$ above and $l$ below is $l=i+1$.
$$
M(X_e\times X_e) \simeq M(X_e) \ds \bigoplus_{l=1}^{e-1} M(X_{e-l,e,e+l})(l^2) \ds M(X_e)(e^2)
$$

This in turn yields a decomposition of the middle-dimensional component of the Chow group of $X_e\times X_e$.
$$
\CH_{e^2}(X_e\times X_e) \simeq \CH_{e^2}(X_e) \ds \bigoplus_{l=1}^{e-1} \CH_{(e-l)(e+l)}(X_{e-l,e,e+l}) \ds \CH_0(X_e)
$$

It now suffices to check the congruence $\mult(\de)\equiv\mult(\de^t) \pmod{2}$ for $\de$ in the image of any of these summands. We treat the first and last summands separately from the rest.

The embedding of the first summand $\CH_{e^2}(X_e)$ is induced by the diagonal morphism $X_e\ra X_e\times X_e$, so the multiplicities of $\de$ and $\de^t$ are equal by symmetry.

For the last summand $\CH_0(X_e)$ we need the following fact.

\begin{proposition}
Any element of $\CH_0(X_e)$ has even degree.
\end{proposition}

\begin{proof}
If $\CH_0(X_e)$ has an element of odd degree, then there exists a field extension $K/F$ of odd degree over which $X_e$ has a rational point.  By \cite[Prop. 1.17]{KMRT98}, $\ind A_K$ divides $e$.  Since the degree of $K$ over $F$ is relatively prime to $\ind A=2e=2^{a+1}$, extension by $K$ does not reduce the index of $A$ \cite[Th. 3.15a]{S99}. Thus $\ind A =\ind A_K$ divides $e$, a contradiction.
\end{proof}

Let the element $\gamma\in\CH_0(X_e)$ have image $\de\in\CH_{e^2}(X_e\times X_e)$.  By the proposition, $\deg (\gamma)$ is even.  For some field $E/F$ over which $X_e$ has a rational point, we set $\bar{X}_e:=(X_e)_E$.  Since $\CH_0(\bar{X}_e)$ is generated by a single element of degree $1$, the image of $\gamma$ in $\CH_0(\bar{X}_e)$ is divisible by $2$.  It follows that $\de\in\CH_{e^2}(\bar{X}_e\times \bar{X}_e)$ is also divisible by $2$ and, since multiplicity does not change under field extension, $\mult (\de)$ is even.  The same argument can be applied to $\de^t$, so $\mult(\de)\equiv 0\equiv \mult(\de^t) \pmod{2}$.

The remaining summands are dealt with by the following proposition.

\begin{proposition}
Let $Fl:=X_{d_1,d_2,\ldots,d_k}(A)$ with $d:=\gcd(e, d_1,d_2,\ldots,d_k)<e$, and let the correspondence $\al:Fl\lt X_e\times X_e$ induce an embedding $$\al_*:\CH_r(Fl)\hookrightarrow \CH_{e^2}(X_e\times X_e).$$  Then for any $\delta$ in the image of $\al_*$, $\mult(\de)\equiv 0\equiv \mult(\de^t) \pmod{2}$.
\end{proposition}

\begin{proof}
Consider the diagram below of fiber products, where we select either of the projections $p_i$ and choose the other morphisms accordingly.

\[
\xymatrix@!C{
                                        & (Fl)_{F(X_e)}  \ar[dr]                       &   \\
(Fl\times X_e)_{F(X_e)} \ar[r] \ar[d] \ar[ur]  & (X_e)_{F(X_e)} \ar[r]  \ar[d]      & \Spec F(X_e) \ar[d]    \\
Fl\times X_e\times X_e \ar[r] \ar[d]    & X_e\times X_e \ar[r]_-{p_2}^-{p_1}    & X_e \\
Fl
}
\]

Taking push-forwards and pull-backs, we get the following diagram which commutes except for the triangle at the bottom.  The push-forward by $p_i$ takes a cycle $\de\in\CH_{e^2}(X_e\times X_e)$ to $\mult(\de)$ if we chose the first projection $p_1$ or to $\mult(\de^t)$ if we chose the second projection $p_2$.

\[
\xymatrix@!C{
      & \CH_0\left((Fl)_{F(X_e)}\right)  \ar[dr]+<-1.5ex,1ex>^-{\deg}                &   \\
\CH_0\left((Fl\times X_e)_{F(X_e)}\right) \ar[r] \ar[ur]   & \CH_0\left((X_e)_{F(X_e)}\right) \ar[r]_-{\deg}  & \Z \ar@{=}[d] \\
\CH_{e^2}\left(Fl\times X_e\times X_e\right) \ar[r] \ar[u]    & \CH_{e^2}\left(X_e\times X_e\right) \ar[r]_-{(\mult)\circ(\trans)}^-{\mult} \ar[u]    & \Z                     \\
\CH_r(Fl) \ar@/_1pc/_-{\al_*}@{.>}[ur] \ar[u]
}
\]

Any $\de\in \im(\al_*)$ also lies in the image of $\CH_{e^2}\left(Fl\times X_e\times X_e\right)$, by the definition of the push-forward.  Chasing through the diagram, one sees that $\mult(\de)$ (and similarly $\mult(\de^t)$) must lie in $\deg \CH_0\left((Fl)_{F(X_e)}\right)$.  We will be done if we can show that no element of $\CH_0\left((Fl)_{F(X_e)}\right)$ has odd degree.

Note that
$$Fl_{F(X_e)} = X_{d_1,d_2,\ldots,d_k}(A)_{F(X_e)} \simeq  X_{d_1,d_2,\ldots,d_k}\left(A_{F(X_e)}\right),$$
where $A_{F(X_e)}$ has index equal to $\gcd(2e,e)=e$ \cite[Th. 2.5]{SV92}.  If some element of $\CH_0\left((Fl)_{F(X_e)}\right)$ has odd degree, then there exists a field extension $K/F(X_e)$ of odd degree over which $(Fl)_{F(X_e)}$ has a rational point.  By Proposition \ref{flagequiv}, $X_d(A_{F(X_e)})$ also has a rational point over $K$.  Thus $\ind A_K$ divides $d<e$, which contradicts $\ind A_{F(X_e)}=e$, since an odd degree extension cannot reduce the index of $A_{F(X_e)}$ \cite[Th. 3.15a]{S99}.
\end{proof}

This completes the proof of the theorem.

\end{document}